\documentclass{amsart}
\usepackage{amssymb}
\usepackage[latin1]{inputenc}
\usepackage[swedish,english]{babel}
\usepackage{ifpdf}
\ifpdf
  \usepackage[pdftex]{graphicx}
  \usepackage{epstopdf}
\else
  \usepackage[dvips]{graphicx}
\fi
\usepackage{subfigure}
\usepackage{sidecap}
\usepackage[pdftitle={Distributed Adaptive FMM3D},
  pdfauthor={Bull, Engblom},
  pdffitwindow=true,
  breaklinks=true,
  colorlinks=true,
  urlcolor=blue,
  linkcolor=red,
  citecolor=red,
  anchorcolor=red]{hyperref}
\usepackage{algorithm}
\usepackage{algpseudocode}
\usepackage[numbers,sort]{natbib}

\numberwithin{equation}{section}
\numberwithin{table}{section}
\numberwithin{figure}{section}

\newcommand{\thetacriterion}{$\theta$-criterion}
\newcommand{\realdom}{\mathbf{R}}
\newcommand{\Ordo}[1]{\mathcal{O}\left(#1\right)}

\begin{document}

\title[Distributed Adaptive FMM in 3D]{Distributed and Adaptive Fast
  Multipole Method In Three Dimensions}

\author[J. Bull]{Jonathan Bull}
\author[S. Engblom]{Stefan Engblom}

\address{Division of Scientific Computing \\
  Department of Information Technology \\
  Uppsala University \\
  SE-751 05 Uppsala, Sweden.}
\urladdr{\url{http://user.it.uu.se/~stefane}}
\email{stefane@it.uu.se}

\thanks{Corresponding author: S. Engblom, telephone +46-18-471 27 54,
  fax +46-18-51 19 25.}





\date{\today}

\begin{abstract}
  We develop a general distributed implementation of an adaptive fast
  multipole method in three space dimensions. We rely on a balanced
  type of adaptive space discretisation which supports a highly
  transparent and fully distributed implementation. A complexity
  analysis indicates favorable scaling properties and numerical
  experiments on up to 512 cores and 1 billion source points verify
  them. The parameters controlling the algorithm are subject to
  in-depth experiments and the performance response to the input
  parameters implies that the overall implementation is well-suited to
  automated tuning.
\end{abstract}

%
%
%

\selectlanguage{english}

\maketitle


\section{Introduction}

The $N$-body problem is one of the most fundamental problems in
computational physics and it has attracted considerable interest from
researchers in numerical algorithms. For the computation of all
pairwise particle interactions, the computational complexity of the
immediate double for-loop algorithm scales as $\Ordo{N^2}$. For large
enough $N$ and high enough tolerance requirements, it is known that
the Fast Multipole Method (FMM) \cite{greengard:87} stands out as an
optimal algorithm of $\Ordo{N}$ complexity, and which also achieves
this optimal performance in practice.  There is a steadily growing
body of research into the use of FMM for the solution of integral
equations and PDEs
\cite{gholami:14,ibeid2018,darve:17,yokota:16,askham:2017}.

FMMs offer a linear complexity in $N$ and enjoy sharp \textit{a
  priori} error estimates, but due to their tree-based nature they are
also notoriously hard to implement, particularly so in three spatial
dimensions \cite{new_AFMM}. With modern multicore and distributed
computer systems, there is an increased interest in designing
effective parallelization strategies
\cite{Chandramowlishwaran:2010b,Cruz2011,
  holm2014,Jinshi2010,Yokota_FMM_manycore,fmmgpu,fmmgpu3,
  exascaleFMM}. Balancing efficiency with implementation transparency
is here an important aspect \cite{ace,pace}.

When source points are non-uniformly distributed it becomes necessary
to adapt the FMM tree according to the local point density such that
the computational effort is balanced across the tree. Mesh adaptivity,
while solved in theory already in the early ages of the FMM
\cite{AFMM}, becomes a major issue when confronted with parallelism
due to the associated complicated memory access patterns. Early
attempts to mitigate this through post-balancing algorithms
\cite{balancedtrees} are less attractive for complexity reasons, and
when data-parallel accelerators made their debut a decade ago, it was
in fact suggested that non-adaptive FMMs offer a better performance
\cite{fmmgpu}.

In this paper we describe the development of a variant of the adaptive
FMM, the balanced tree FMM \cite{engblom:2011, goude:2013, holm2014},
in three dimensions and using distributed parallelisation on very
large modern computer systems. Rather than adapting the number of
levels locally as in the classic adaptive FMM, our balanced tree FMM
maintains a fixed number of levels and splits boxes at median planes
such that the number of points in each subtree is balanced at every
level. Otherwise known as orthogonal recursive bisection (ORB)
\cite{warren1992,singh1993,yokota-petascale}, it is guaranteed to
produce a balanced tree \cite{yokota-petascale}.

The implementation of the balanced tree FMM with distributed
parallelism is relatively straightforward and is shown to scale well
up to 512 processes at a respectable absolute efficiency.
Furthermore, it exposes various parameters to the user enabling the
performance to be readily optimised for a particular simulation. We
show that the performance response is convex in the parameter space
making the balanced tree FMM particularly well-suited to automatic
tuning in an online computing context.

Since the FMM has been judged to be one of the top 10 most important
algorithms of the 20th century \cite{Top10algorithms}, it is our hope
that insights obtained here is of general value. Our implementation of
the parallel balanced tree FMM is therefore freely available as
\texttt{daFMM3D} (distributed adaptive FMM in 3D) C/C++/Matlab code at
\url{www.stenglib.org}.

The structure of the paper is as follows: in \S\ref{sec:algorithm} we
describe the balanced tree FMM in three dimensions including the
multipole acceptance criterion, adaptive box splitting and
computational complexity estimates. The corresponding distributed
algorithm is also presented with a description of parallel data
structures and communication complexity. In \S\ref{sec:examples}
convergence test results, parameter response test results and and
strong scaling results on up to 512 processes are reported.  The
algorithm's adaptive response is also tested in the more challenging
case of a spiral galaxy of 1 billion sources with a multiscale
structure.


\section{3D balanced tree FMM in a distributed environment}
\label{sec:algorithm}

Fast multipole methods evaluate pairwise interactions of the type
\begin{align}
  \label{eq:paireval}
  \Phi(x_{i}) &= \sum_{j = 1, j \not = i}^{N} G(x_{i},x_{j}),
  \quad x_{i} \in \realdom^{D}, \quad i = 1 \ldots N,
\end{align}
where $D = 3$ in this paper, and where $\{x_{j}\}$ are a set of $N$
\emph{sources} (points) in a force field governed by the kernel
$G$. The FMM is a tree-based algorithm that produces a continuous
representation of the field such that it can be evaluated to within
the tolerance anywhere inside a bounding box enclosing the sources.

\subsection{Adaptivity and complexity}

The sources are initially placed in a single bounding box at level
$l = 0$ of the tree. Splitting operations produce child boxes on
levels $l > 0$, with the number of points per box becoming
successively smaller with each level. Each box is given an outgoing
\textit{multipole} expansion and an incoming \textit{local} expansion.
The multipole expansion is valid far away from the box and is the
expansion of all the points in the box about the box centre. In 3D it
is the spherical harmonics expansion \cite{greengardshortcourse}. The
local expansion is valid within the box and is the expansion of
far-away sources about the box centre. Both the multipole and local
expansions are of order $Q$ which determines the overall precision.

\subsubsection{The \thetacriterion}

We consider the \textit{balanced tree} version of FMM as described in
\cite{engblom:2011}.  As with the original FMM, boxes on the same
level are either unconnected or strongly/weakly connected depending on
the separation of their centres.  A pair of boxes with radii $r_1$ and
$r_2$ and separation $d$ are deemed to be weakly connected if the
\textit{$theta$ criterion} is satisfied:
\begin{equation}
  \label{eq:thetacrit}
  R + \theta r \leq \theta d,
\end{equation}
where $R$ = max$(r_1,r_2)$, $r$ = min$(r_1,r_2)$ and $\theta \in
(0,1)$ is a parameter.  If the criterion is not satisfied then the
boxes are strongly connected.  The choice of $\theta$ affects the
\emph{connectivity pattern}, i.e., the stencil of strongly and weakly
connected boxes, and thus the relative amount of work spent on
different parts of the algorithm. In practice it has been found that
$\theta = 0.5$ leads to a well balanced algorithm in 2D for various
source densities \cite{goude:2013}, however, we will reevaluate the
optimal value for different point distributions in 3D in this
paper. Boxes inherit connections from their parents: the children of
strongly connected parents are allowed to be strongly or weakly
connected to each other, as judged by \eqref{eq:thetacrit}, whereas
the children of weakly connected parents are not connected.

\subsubsection{Adaptivity}

A balanced tree, or a \emph{pyramid}, is one in which all boxes at
level $l$ are either subdivided further into level $l+1$, or form the
\emph{leafs} of the tree, i.e., the boxes at the finest
level. Traditional level-adaptive FMMs \cite{new_AFMM,AFMM} work by
adjusting the number of levels locally, resulting in equal-sized boxes
and an unbalanced tree. The balanced tree FMM \cite{engblom:2011}, by
contrast, keeps the number of levels fixed but each box is adaptively
split along median planes such that the child boxes contain about an
equal number of points, implying a pyramid tree. This results in
different-sized boxes and a variable connectivity pattern, but crucially
avoids the cross-level communication required by level-adaptive FMM.

For additional control over the tree we introduce an adaptivity
parameter, $\eta \in [0,1]$.  $\eta = 0$ enforces splitting at
geometric mean planes resulting in a non-adaptive or fixed tree where
all children are the same size and there are a uniform number of
connections between boxes, thus balancing the cost of the
Multipole-to-Local (M2L) shift. $\eta = 1$ instead splits boxes at
median planes such that all children contain about an equal number of
sources, thus balancing the cost of the Particle-to-Particle (P2P)
shift.  One may set $\eta$ in between 0 and 1 to split at a weighted
average of the geometric mean and the source point median.  This may
be useful for tuning the performance in individual simulations.  A
similar approach splitting at the $n$th element has been used to
balance loads between arbitrary numbers of parallel processes
\cite{yokota-petascale}.

\subsubsection{Computational Complexity}

The cost of the FMM is dominated by the M2L and P2P shifts.  We
evaluate the complexity of these two operations in the 3D balanced
tree FMM with $\eta=1$ and $\theta$ arbitrary.  It is assumed that the
mesh is asymptotically regular so that $R \approx r$ in
\eqref{eq:thetacrit} and that there are $M$ boxes at the finest level.

Each box is weakly connected to boxes within a spherical shell whose
inner surface is defined by \eqref{eq:thetacrit},
$d = (1+\theta)/\theta \times r$.  The outer surface is defined by
strong connections at the parents' level,
$d_{\mbox{{\tiny parent}}} = (1+\theta)/\theta \times r_{\mbox{{\tiny
      parent}}} = 2d$. This shell contains approximately $d^3 M$
boxes. By approximating the separation $r \sim M^{-1/3}$, we arrive at
an estimate of $\Ordo{\theta^{-3}}$ weak connections per box. M2L in
its basic form scales as $\Ordo{Q^4}$ \cite{greengardshortcourse}.  We
use the rotational translation form which rotates the $z$-axis of the
spherical harmonic expansion to the vector from one box to another.
This form requires $\Ordo{Q^3}$ operations per pair so the complexity
of the M2L shift for all $M$ boxes is $\Ordo{Q^3 \theta^{-3} M}$.  The
exponential plane wave expansion form scales as $\Ordo{Q^2}$
\cite{greengardshortcourse}, but it is not practical to use in the
balanced tree FMM because the geometric relationships are not fixed
\textit{a priori}.

Each box is strongly connected to $\frac{4}{3}\pi d^3 M$ boxes in a
sphere, each containing $N/M$ sources. The P2P evaluation requires
$\Ordo{N^2 M^{-2}}$ operations per pair of boxes so the total cost of
P2P is $\Ordo{N^2 \theta^{-3} M^{-1}}$.  If we always set the number
of levels such that $N/M$ is constant this becomes
$\Ordo{\theta^{-3} N}$. To balance the M2L and P2P operations we
require that the number of source points per leaf is related to the
order of the multipole expansion: $(N/M)^2 \propto Q^3$.

In practice, the assumption of asymptotic regularity may not be
satisfied, particularly for low numbers of levels and highly
inhomogeneous source distributions, leading to different scaling than
the above estimates.  However, with the free parameters $L$ (which
determines the ratio $N/M$), $\theta$ and $\eta$ it is possible to
tune the balanced tree FMM for a specific source distribution.

\subsection{Distributed adaptive FMM}

This section describes the parallelisation of the balanced tree FMM
for $N$ points with $L+1$ levels $l \in \{0,1,2,\ldots,L\}$ in a
distributed memory framework on $P$ processes. The balanced tree FMM
requires the construction of a connectivity matrix since geometric
relationships between boxes are not fixed \textit{a priori}. In
parallel, this connectivity matrix has to be partially shared between
processes.

\subsubsection{Data Structures}

In the FMM octree $T$ we refer to the box on level 0, denoted
$T^{(0)}$, as the global root.  It is split into $P$ boxes on level 1,
the local roots $T_p^{(1)}$.  Each process $p$ then constructs its own
$L$-level local octree $T_p$ with local root $T_p^{(1)}$ and leaf
level $T_p^{(L)}$. The $j$th box on level $l$ is denoted $T_{p,j}^{l}$.

\begin{figure}[hbt]
  \includegraphics[width=0.8\textwidth]{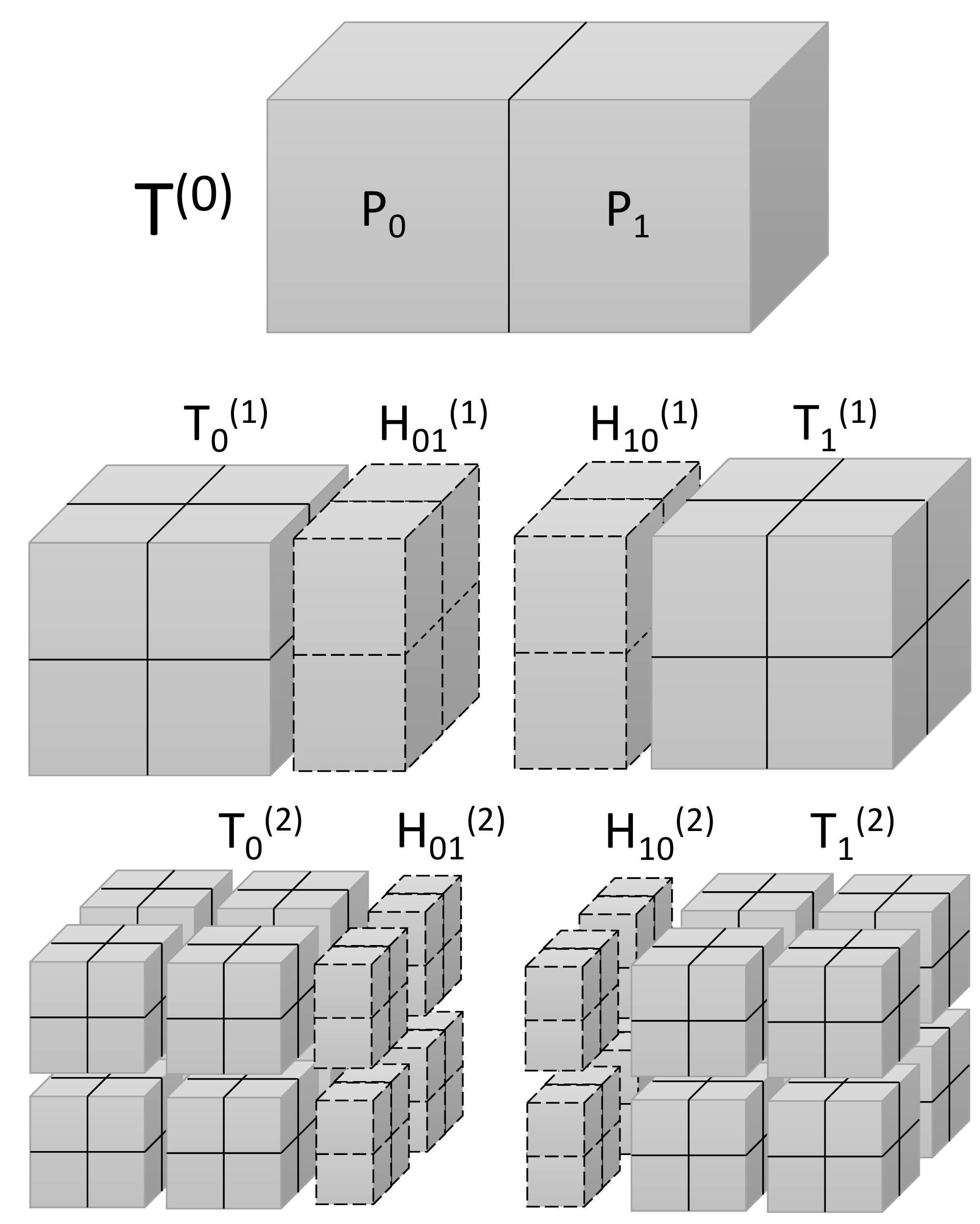}
  \caption{Schematic of a 2-process distributed tree with 3
    levels. Halos are shown as boxes with dashed outlines. Subscripts
    denote parallel processes and superscripts denote tree level.}
  \label{fig:tree}
\end{figure}

Since the connectivity pattern is not fixed in the balanced tree FMM,
supporting data structures are required. On levels $l=1:L$ each
process $p$ constructs a local compressed column storage (CCS) format
sparse connectivity matrix $C_p^{(l)}$ encoding same-level connections
within $T_p$, as well as $P-1$ (possibly empty) halo matrices
$H_{pq}^{(l)}$ encoding same-level connections to $T_q, q \neq p$.
Together, all $C^{(l)}$ and $H^{(l)}$ matrices form a global
connectivity matrix on level $l$.  The halos are symmetric
($H_{qp}^{(l)} = (H_{pq}^{(l)})^T$) so only the lower-triangular part
of the global matrix needs to be constructed, however to reduce
communication each process constructs all $P-1$ halos itself.  Figure
\ref{fig:tree} illustrates a distributed FMM octree.

Each halo $H_{pq}^{(l)}$ contains a sparse matrix in index-format that
stores all the pairs of strongly and weakly connected boxes owned by
processes $p$ and $q$. The index-sparse format utilises two lists,
\texttt{ibox} and \texttt{jbox}, containing the indices of boxes in
$p$ and $q$, respectively. In the FMM, boxes can only interact if
their parents are strongly connected, hence
$H_{pq}^{(l)}$.\texttt{ibox} is 64 times larger than
$H_{pq}^{(l-1)}$.\texttt{ibox}; likewise for
$H_{pq}^{(l)}$.\texttt{jbox}. The index-sparse format implies that a
box appears $k$ times if it has $k$ connections.  The format has been
chosen since it compresses better than the CCS format when only a
small subset of the boxes are connected. This condition is true for
large $l$ where only a very thin layer of boxes on the subdomain
surface are connected across the process boundary.

\subsubsection{Distributed algorithm}

Algorithm \ref{algorithm:FMM} describes the whole distributed balanced
tree FMM for $P$ processes and an $L+1$ level octree.  
Construction of the balanced tree, connectivity matrix and halo
matrix (Tree subroutine) are described in Algorithm \ref{algorithm:FMMh}
in the Appendix. The downward pass of the balanced tree FMM (M2L, L2L
subroutines) is described in Algorithm \ref{algorithm:FMMd}. 
The direct evaluation (P2P, L2P subroutines) is described in Algorithm
\ref{algorithm:FMMe}. It is assumed that data partitioning onto $P$
processes is done in a separate preprocessing step, mimicking the
intended use as a plugin solver for other scientific applications.

\begin{algorithm}
  \caption{Distributed balanced tree FMM.}
  \label{algorithm:FMM}
  \begin{algorithmic}[1]
    \ForAll{$p \in P$}
    \Statex \textit{Alloc}:
    \State allocate arrays for tree $T_p$ on all levels and for points on level $L$
    \Statex \textit{Tree} (construction of $T_p$, $C_p$, $H_p$):
    \State Algorithm \ref{algorithm:FMMh}
    \Statex \textit{P2M}:
    \State Compute multipole expansions on level $L$ from point data
    \Statex \textit{M2M} (Upward pass of tree $T_p$):
    \For{$l=L:2$}
    \State compute multipole expansions in $T_p^{(l-1)}$ from children in $T_p^{(l)}$
    \EndFor
    \Statex \textit{M2L, L2L} (Downward pass of tree $T_p$):
    \State Algorithm \ref{algorithm:FMMd}
    \Statex \textit{L2P, P2P} (Direct valuation of potentials on leaf level $L$):
    \State Algorithm \ref{algorithm:FMMe}
    \EndFor
  \end{algorithmic}
\end{algorithm}

\subsubsection{Computational Complexity}

Besides the M2L and P2P shifts, which dominate the cost of FMM, there
is a cost associated with halo construction in the distributed
balanced tree FMM.  Halo construction, detailed in Algorithm
\ref{algorithm:FMMh}, has computational complexity of $\Ordo{N}$.

\subsubsection{Communication Complexity}

Communication is dominated by the P2P interaction amongst $N$ sources.
Let the domain with $\Ordo{1}$ volume be divided equally into $P$
subdomains and assume that $P$ is sufficiently large that relatively
few processes lie on the domain boundary. The communication complexity
is proportional to the number of sources lying close to each
subdomain's boundaries with other subdomains. Each subdomain has a
surface area of $\Ordo{P^{-2/3}}$ and there are $\Ordo{N^{2/3}}$ 
sources per unit area.  This leads to the
estimate of $\Ordo{(N/P)^{2/3}}$ for the communication complexity as
calculated in \cite{yokota:2014a}. There is also a contribution of
$\Ordo{(M/P)^{2/3}}$ for the M2L shift between $M$ boxes but since we
let $N/M$ be constant this is included in the above
estimate. Similarly the halo requires communication of the positions
and dimensions of $M$ boxes which can be included in the above
estimate.

\subsubsection{Summary}

The distributed balanced tree FMM is designed to allow easy
tuning/optimisation for a particular distribution of sources by
exposing the parameters of relevance to the algorithm. It is
relatively straightforward to implement and memory access is simpler
than in the classic level-adaptive FMM.  The additional cost of halo
construction and communication is negligible.


\section{Computational examples}
\label{sec:examples}

The distributed balanced tree FMM has been implemented in a C/C++/MPI
code named \texttt{daFMM3D} available at \url{www.stenglib.org}. In
serial it can be run via a Matlab interface. The code solves $N$-body
problems subject to forces that follow an inverse-square law such as
gravitation, although other kernels can be implemented analogously.

\subsection{Convergence}

\begin{figure}[hbt]
  \includegraphics[width=0.8\textwidth]{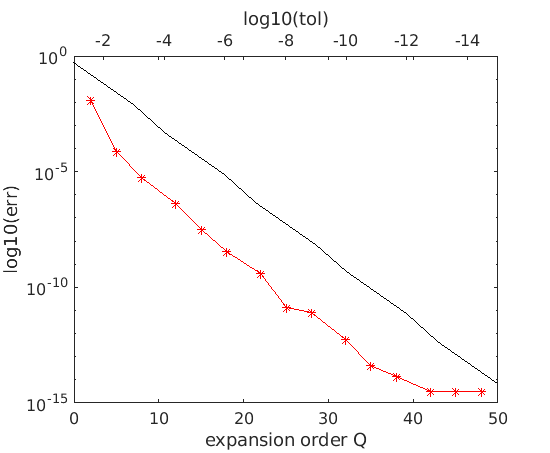}
  \caption{Convergence of relative error in FMM solution. The upper
    $x$-axis shows the tolerance and the lower $x$-axis the polynomial
    order of the expansions.}
  \label{fig:convergence}
\end{figure}

Engblom (2011) \cite{engblom:2011} derived that the relative error of
the balanced tree FMM solution $\phi$ for $N$ positive potentials
under the theta criterion is bounded by a factor $C
\theta^{p+1}/(1-\theta)^2$, where $C$ is a constant.  We verify that
this bound is satisfied by the present implementation for a set of $N
= 1000$ random points at $N$ different random evaluation points in
3D. A reference solution $\phi_0$ is computed by direct evaluation and
the FMM solutions $\phi_k$ are computed with tolerances of $10^{-k}$
for $k = 1,2,\ldots,15$. Figure~\ref{fig:convergence} shows the
relative error $e_k = | (\phi_0-\phi_k)/\phi_0 |$ together with the
theoretical bound.

\subsection{Parameter response}

The sensitivity of the balanced tree FMM to the free parameters
$\theta$, $\eta$ and $L$ was investigated for four different
distributions of 100,000 points on a single process. In this study
the error tolerance was fixed to {\tt TOL}$= 10^{-6}$.
The four distributions were (a) randomly distributed points in the unit
cube, (b) spherical Gaussian distribution, (c) spherical shell and (d)
a helix.  These distributions were chosen with the intent to
`stress-test' the adaptive response of the algorithm. The three
independent parameters were varied in the ranges
$0 \leq \theta \leq 0.8$, $1 \leq L \leq 4$, $0 \leq \eta \leq 1$.
While varying one parameter the remaining two were set to default values
($\theta=0.5$, $\eta=0.5$ and $L=3$) so although this is not a full
sweep of parameter space, it gives a reasonably good understanding of
the performance response.

\begin{figure}[hbt]
  \includegraphics[width=0.8\textwidth]{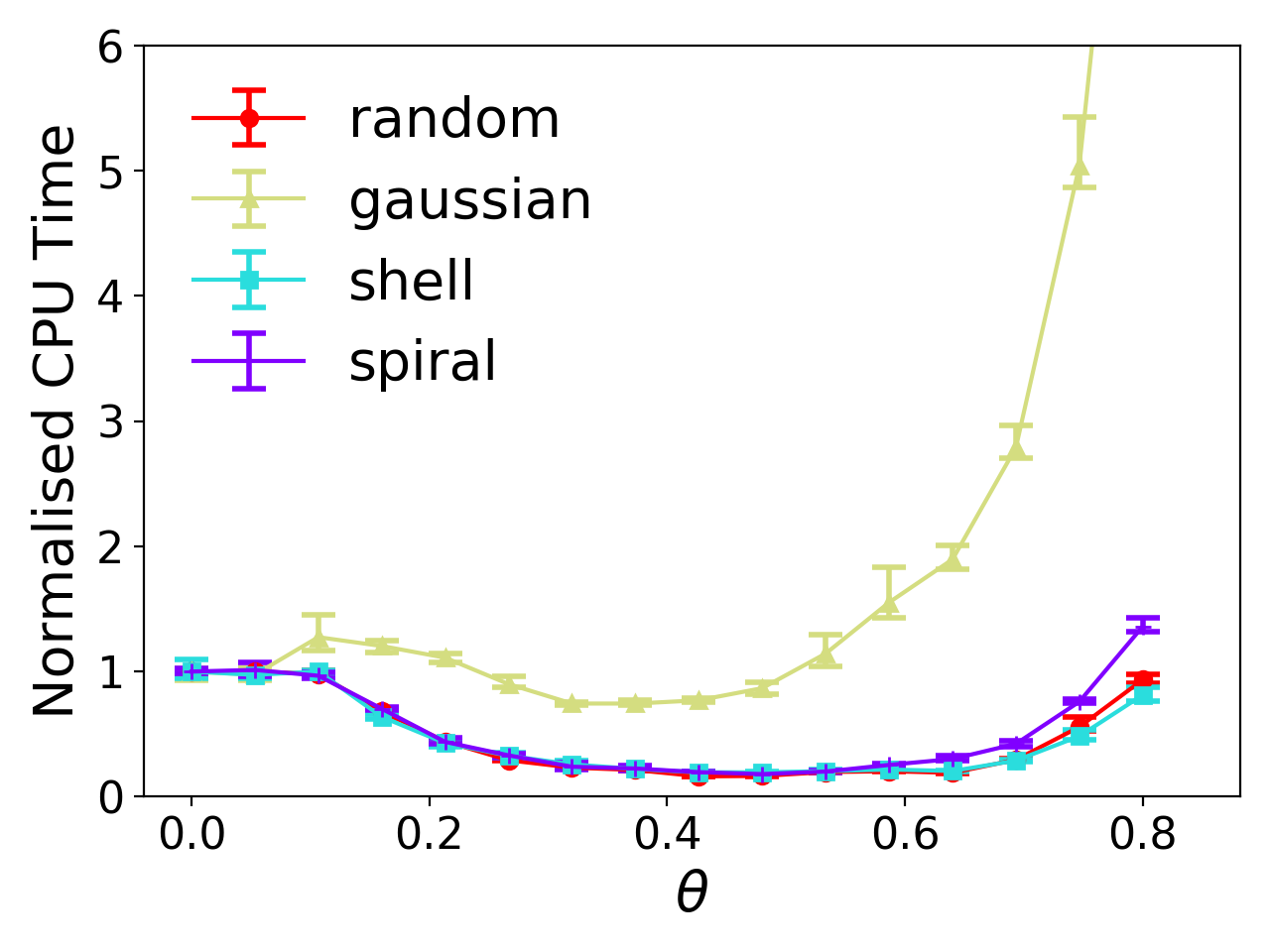}
  \caption{Sensitivity of total time to $\theta$ in the random, Gaussian,
  spiral and shell distributions of 100000 points.}
  \label{fig:thetaopt}
\end{figure}

\begin{figure}[hbt]
  \includegraphics[width=0.8\textwidth]{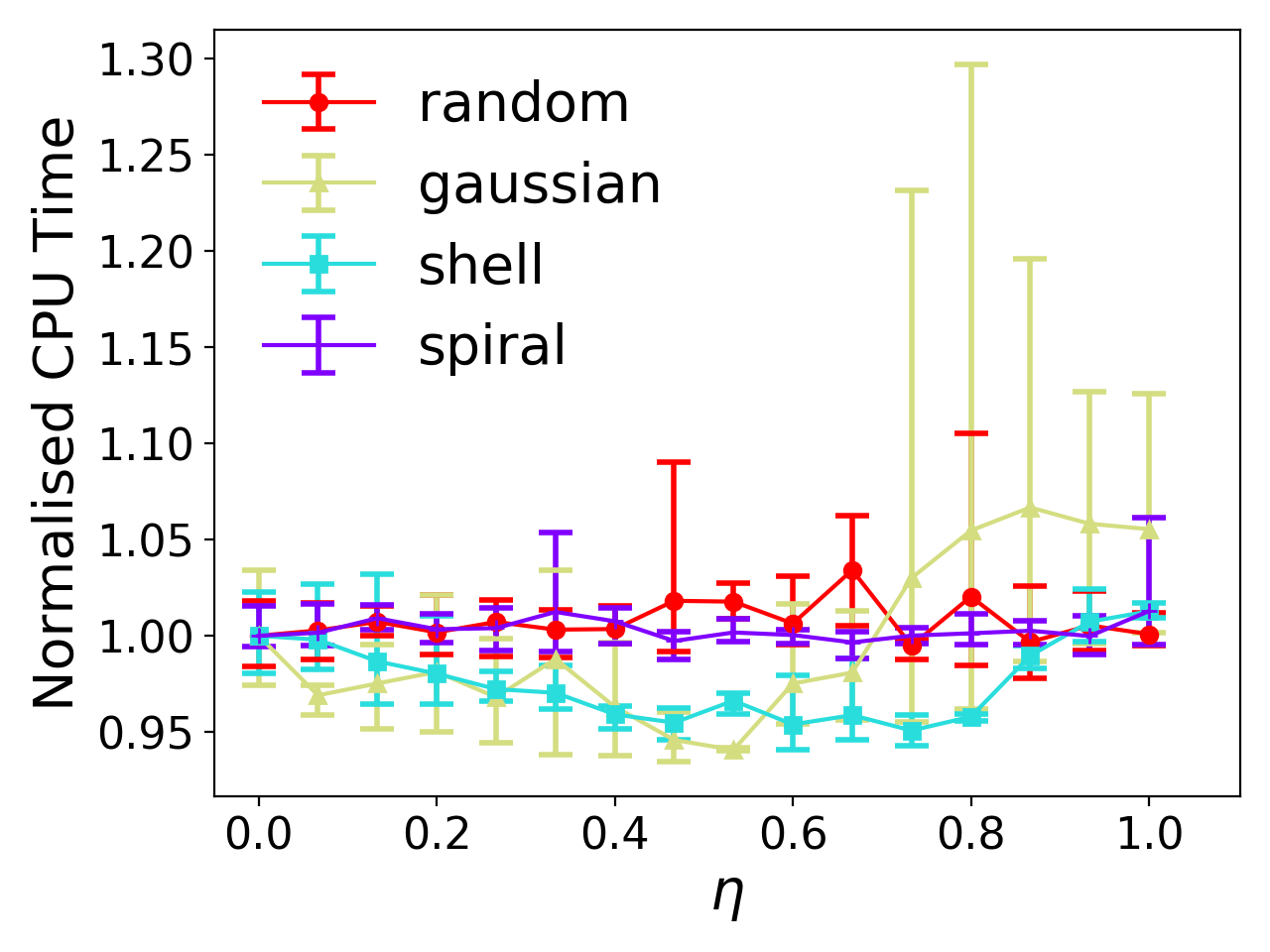}
  \caption{Sensitivity of total time to $\eta$ in the random, Gaussian,
  spiral and shell distributions of 100000 points.}
  \label{fig:etaopt}
\end{figure}

\begin{figure}[hbt]
  \includegraphics[width=0.8\textwidth]{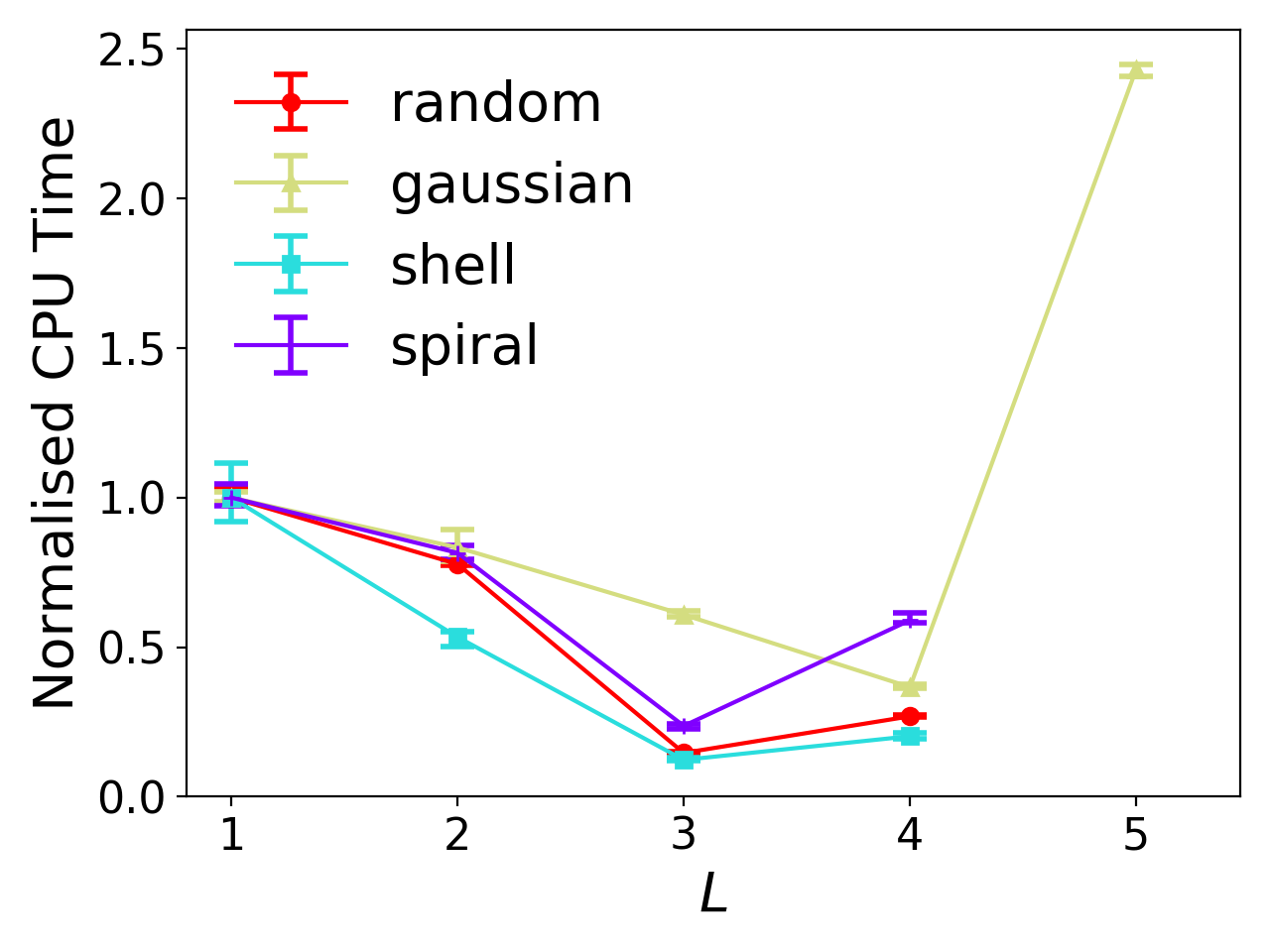}
  \caption{Sensitivity of total time to $L$ in the random, Gaussian,
  spiral and shell distributions of 100000 points. 
  Note 5 values for Gaussian and 4 for the other cases.}
  \label{fig:Lopt}
\end{figure}

The code was executed 4 times for each combination of parameters.
Figures~\ref{fig:thetaopt}, \ref{fig:etaopt} and \ref{fig:Lopt} show
the normalised mean time curves against $\theta$, $\eta$ and $L$
respectively, with error bars denoting the minimum and maximum
times. Times were normalised by the mean times at $\theta=0.5$,
$\eta=0$ and $L=1$. Table~\ref{tab:paramtest} lists the optimal
parameter values identified from the graphs.

All cases are highly sensitive to $\theta$. Run times for the random,
shell and spiral cases are minimised in the range $0.4 < \theta < 0.6$,
consistent with a previous result in 2D \cite{engblom:2011}. 
The Gaussian case has minimum run time at $0.3 < \theta < 0.5$ and it
is less sensitive to low $\theta$ than the other cases.
Interestingly, none of the cases display a strong sensitivity to
$\eta$. The Gaussian and shell cases have distinct minima (at
$\eta \approx 0.5$ and $0.4 < \eta < 0.8$ respectively) but the times
are only 5\% less than at $\eta=0$. 
The Gaussian case exhibits more variation in run times at large values
of $\eta$ than the other 3 cases.
The optimal number of levels is 3 for the random, spiral and shell
cases, corresponding to a maximum of 195 points per leaf box.
Four levels is optimal for the Gaussian, about 24 points per leaf box. 

The optimal values are a compromise between the M2L and P2P subroutines. M2L time increases with $L$ and $\theta$, whereas P2P time decreases
with $L$ and $\theta$. Furthermore, the optimal values of $L$ and
$\theta$ are co-dependent: more levels increases the M2L time which
biases the optimal $\theta$ towards lower values (reduced M2L time)
and vice versa.

These results demonstrate that for a given distribution, it is
generally not possible to determine the optimal parameters \textit{a
  priori}.  However, the clear convex response of total time to
individually varying $L$ and $\theta$ motivates the use of autotuning
optimisation techniques where the parameters are adjusted in steps as
the underlying simulation proceeds. $L$ has the most drastic effect
and could possibly be chosen \textit{a priori} to result in a desired
mean number of points per leaf (suggested range 50-200), leaving
$\theta$ and $\eta$ free to vary for fine tuning.  The four tested
distributions displayed little sensitivity to $\eta$, but this may be
due to their limited size and complexity. In larger, more complex
distributions we expect that $\eta$ may play a larger role, as seen
with the spiral galaxy in \S\ref{subsec:galaxy}.

\begin{table}
  \begin{tabular}{rlll}
    \hline
    distribution & $\theta$ opt. & $\eta$ opt. & $L$ opt. \\
    \hline
    random & 0.45 & - & 3 \\
    Gaussian & 0.4 & 0.5 & 4 \\
    shell & 0.5 & 0.4-0.8 & 3 \\
    helix & 0.5 & - & 3 \\
    \hline
  \end{tabular}
  \caption{Optimal parameter values in the random, Gaussian, spiral and
  shell distributions of 100000 points. A dash (-) means no optimum value
  could be identified (flat response curve).}
  \label{tab:paramtest}
\end{table}

\subsection{Scaling}

Strong and weak scaling is assessed for the \texttt{daFMM3D} code
in non-adaptive mode ($\eta$=0, $\theta$=0.5) on a set of uniformly
distributed
points in the unit cube.  They are partitioned into cubic numbers of
processes ($1^3, 2^3, \ldots$) so that the connectivity patterns are
as close to identical as possible.  To verify that balanced-tree
adaptivity has no impact when the points are uniformly distributed the
weak scaling tests were also run with $\eta=1$.  The total times and
scaling efficiencies were near-identical to the non-adaptive
algorithm.  Every data point in the results is an ensemble average
from a minimum of four separate runs.

Scaling tests were run on the Rackham cluster, part of the Uppmax
computing service at Uppsala University.  Rackham consists of 334
nodes each with two 10-core Xeon E5-2630 V4 processors running at 2.2
GHz (turbo 3.1 GHz) and 128 GB memory per node.  All nodes are
interconnected with a 2:1 oversubscribed FDR (56 GB/s) Infiniband
fabric.  When the requested number of processes fitted onto one node,
that node's entire resources were reserved for that job, eliminating
timing inconsistencies due to sharing the node with other jobs.

\subsubsection{Strong Scaling Results}

\begin{figure}[hbt]
  \subfigure{
    \includegraphics[width=0.7\textwidth]{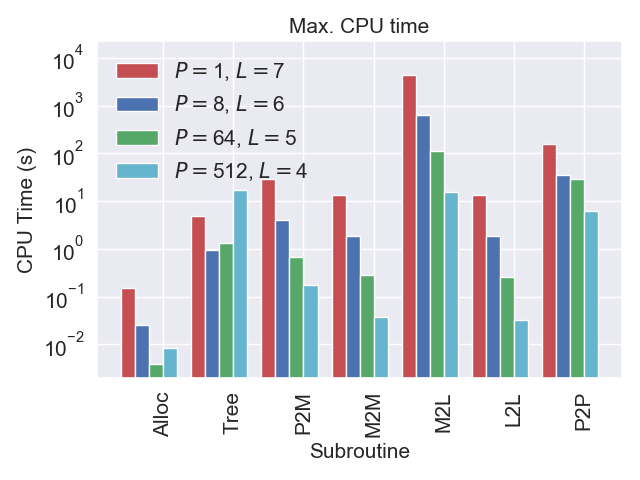}}\\
  \subfigure{
    \includegraphics[width=0.7\textwidth]{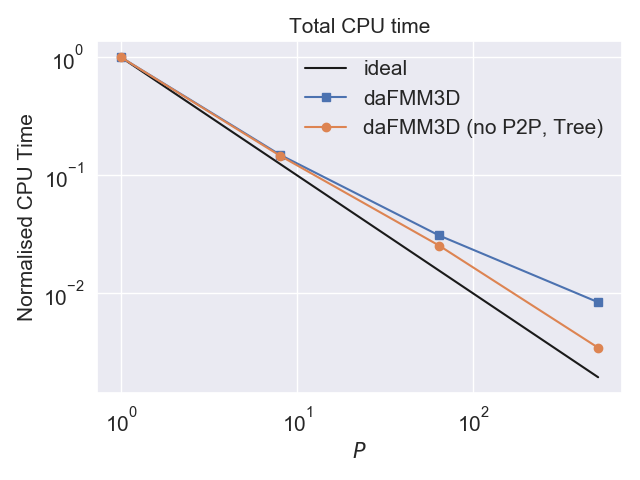}}
  \caption{Strong scaling, 10M points, 1-512 processes, constant 38
    points/leaf box. Top: time per subroutine, bottom: total time.}
  \label{fig:strong}
\end{figure}

In the strong scaling tests $N$=10M points were distributed over 1,
8, 64 and 512 processes with 7, 6, 5 and 4 levels in the global tree
respectively. To ensure a fair comparison the number of levels was
reduced as $P$ increased so that the number of points per leaf box
$N/M$ was kept approximately constant ($N/M \approx
38$). Figure~\ref{fig:strong} shows the CPU time of the different
subroutines as well as the total time. The strong scaling efficiency
was 23.2\% in these tests, and was mainly limited by the non-ideal
scaling of Tree and P2P. Excluding
these, the scaling efficiency of 56.2\% was achieved. M2L was the
dominant subroutine in these tests and it achieved an efficiency of
56.4\%.

\subsubsection{Weak Scaling Results}

\begin{figure}[hbt]
  \subfigure[Max time for all subroutines]{
    \includegraphics[width=\textwidth]{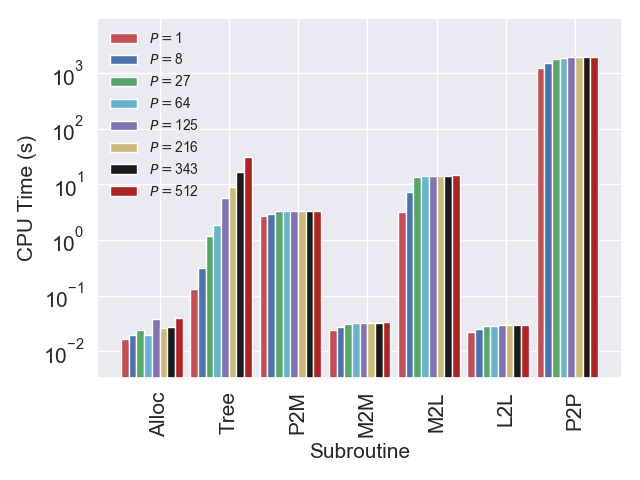}}\\
  \subfigure[Normalised max time for M2L]{
    \includegraphics[width=0.48\textwidth]{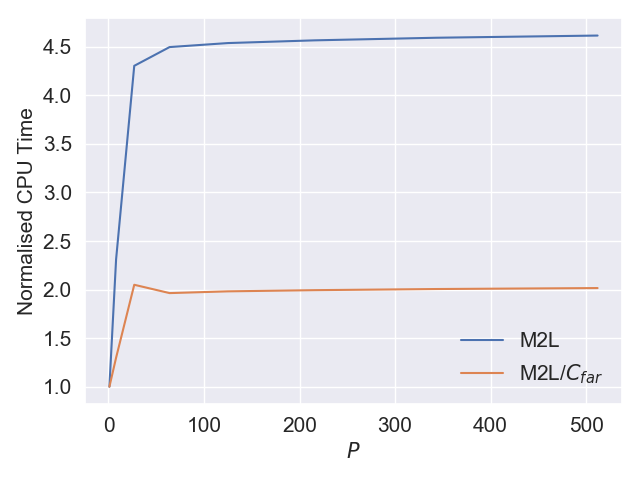}}
  \subfigure[Normalised max time for P2P]{
    \includegraphics[width=0.48\textwidth]{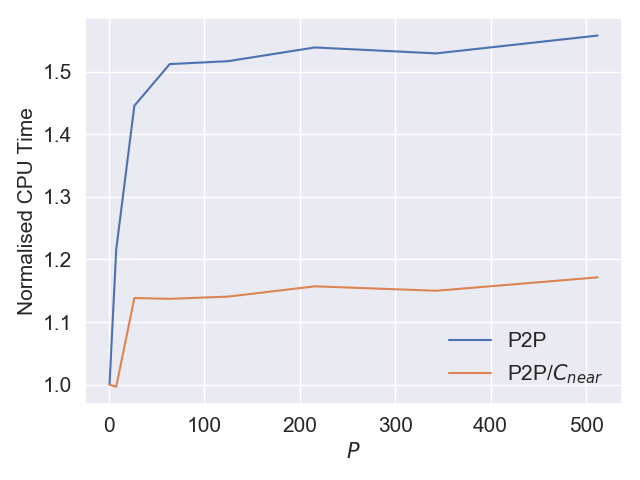}}
  \caption{Weak scaling, 1M points/process, 1-512 processes, 4 levels.}
  \label{fig:weakl4}
\end{figure}

\begin{figure}[hbt]
  \subfigure[Max time for all subroutines]{
    \includegraphics[width=\textwidth]{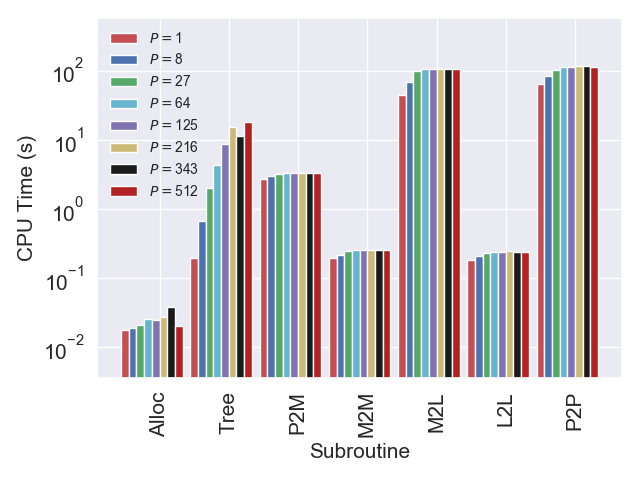}}\\
  \subfigure[Normalised max time for M2L]{
    \includegraphics[width=0.48\textwidth]{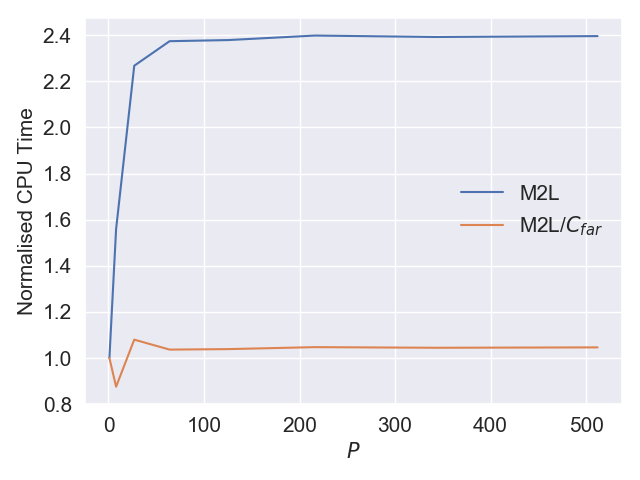}}
  \subfigure[Normalised max time for P2P]{
    \includegraphics[width=0.48\textwidth]{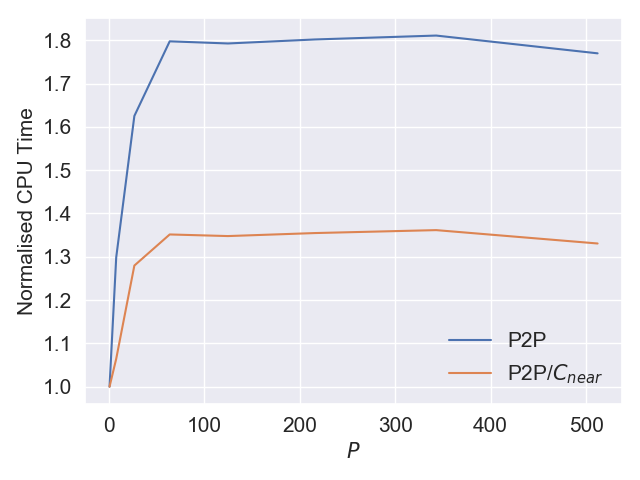}}
  \caption{Weak scaling, 1M points/process, 1-512 processes, 5 levels.}
  \label{fig:weakl5}
\end{figure}

\begin{figure}[hbt]
  \subfigure[Max time for all subroutines]{
    \includegraphics[width=\textwidth]{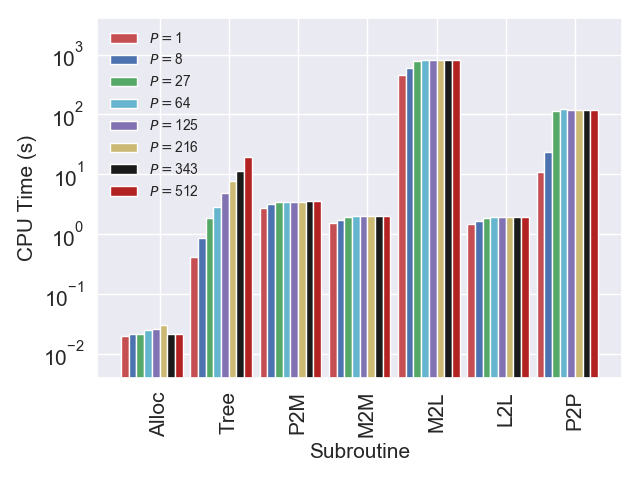}}\\
  \subfigure[Normalised max time for M2L]{
    \includegraphics[width=0.48\textwidth]{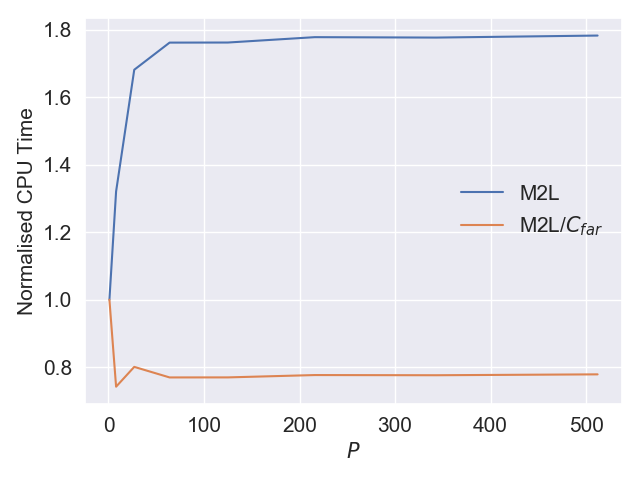}}
  \subfigure[Normalised max time for P2P]{
    \includegraphics[width=0.48\textwidth]{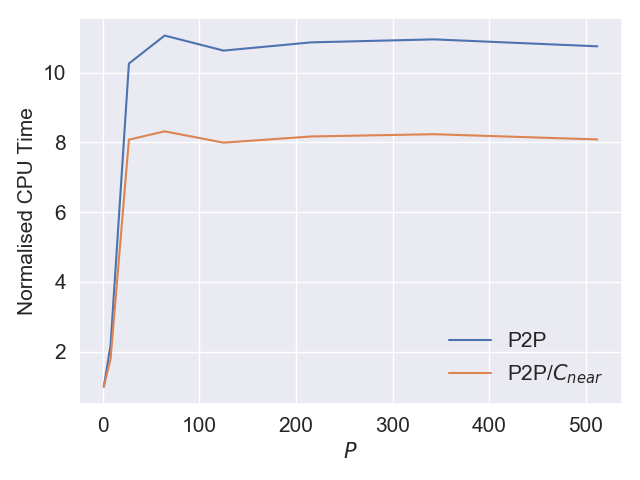}}
  \caption{Weak scaling, 1M points/process, 1-512 processes, 6 levels}
  \label{fig:weakl6}
\end{figure}

In the weak scaling tests 1M points were assigned to each process and
up to 512 processes were used. In Figures~\ref{fig:weakl4},
\ref{fig:weakl5}, and \ref{fig:weakl6}, (a) the maximum CPU times of each
subroutine, (b) the normalised M2L and (c) normalised P2P times
with 4, 5 and 6
levels, respectively, are reported. The total weak scaling efficiency
of 62.9\%, 46.2\% and 49.4\% is achieved with 4, 5 and 6 levels. The
total time is strongly dependent on the number of levels: 5 levels
(244 points/leaf) is about 4 times faster than 4 levels (1953
points/leaf) or 6 levels (31 points/leaf).

All subroutines except Tree achieve relatively good weak scaling
with flat maximum times for $P>$27. The increase in times
from $P$=1 to $P$=27 is anticipated due to connectivity changes as
discussed below. Poor scaling of the Tree subroutine 
is due to latency caused by waiting for non-blocking MPI receive commands
to complete in process number order (Algorithm \ref{algorithm:FMMh} 
line 22). Tree makes a negligible contribution to total time within the
range of $P$ tested but it will become significant at very large $P$;
this is easily mitigated by some additional software engineering.
The M2L and P2P subroutines were implemented differently, with 
non-blocking receives being processed in the order in which they were completed to reduce latency 
(see Algorithm \ref{algorithm:FMMd} lines 7-15 and 
Algorithm \ref{algorithm:FMMe} lines 10-17). 

The number of levels influences
the scaling of M2L and P2P differently: the best weak scaling of M2L
is with 6 levels but the best scaling of P2P is with 4 levels.
Conversely the shortest execution time for M2L is with 4 levels and
for P2P with 6 levels. 5 levels offers the best balance of scaling
efficiency and execution time for both subroutines and the lowest
total execution time.

It is not possible to obtain perfect scaling of the FMM algorithm due
to connectivity changing with the number of processes: as $P$
increases a larger fraction of the subdomains do not have a boundary
on the domain outer boundary so the average number of connections
(especially far connections) to other processes increases. The
largest growth in connectivity is from $P$=1 to $P$=27 and
thereafter the increase in connectivity is marginal.
Denoting the total number of near connections on the leaf level as
$N_{near}$ we define a near connectivity factor for $P$ processes:
$C_{near}=N_{near}/N_{near,P=1}$.
Likewise the far connectivity factor is
$C_{far}=N_{far}/N_{far,P=1}$ where $N_{far}$ is the total
number of far connections (all levels).
Table~\ref{tab:connectivity_weak} lists $C_{near}$ and $C_{far}$
for the 4-level weak scaling test up to $P$=64. To avoid burdensome
calculations we assume that the 4-level, 64-process factors are good
estimates of the factors that would be obtained with more levels and
more processes.

In Figures \ref{fig:weakl4}-\ref{fig:weakl6} the orange lines show (b)
the M2L times divided by $C_{far}$ and (c) the P2P times divided by
$C_{near}$ for 4, 5 and 6 levels respectively.
Table~\ref{tab:weakscaleff} lists the scaling efficiencies of M2L, P2P
and the total algorithm without and with adjustment by $C_{near}$ and
$C_{far}$. The adjusted efficiencies of the total algorithm are
83.5\%, 72.7\% and 100.6\% for 4, 5 and 6 levels respectively. The
adjusted efficiencies of M2L for $L>4$ are greater than 100\%,
suggesting that some of the extra connectivity burden at $P>1$ may be
hidden by efficiency gains obtained through good parallel execution
(non-blocking calls and hiding latency with local computation). On the
other hand, the adjusted P2P efficiencies are still low, showing that
there is room for improvement in the parallelisation of this
subroutine.

\begin{table}
  \begin{tabular}{rll}
    \hline
    $P$ & $C_{near}$ & $C_{far}$ \\
    \hline
    1 & 1 & 1 \\
    8 & 1.22 & 1.78 \\
    27 & 1.27 & 2.10 \\
    64 & 1.33 & 2.29 \\
    \hline
  \end{tabular}
  \caption{Growth of total number of near and far connections with $P$
  relative to $P$=1 for 4-level weak scaling test setup.}
  \label{tab:connectivity_weak}
\end{table}

\begin{table}
  \begin{tabular}{rlll}
    \hline
    levels & M2L (adjusted) & P2P (adjusted) & Total (adjusted) \\
    \hline
    4 & 21.7 (49.6) & 64.2 (85.4) & 62.9 (83.5) \\
    5 & 41.7 (95.6) & 56.5 (75.1) & 46.2 (72.7) \\
    6 & 56.1 (128.4) & 9.3 (12.4) & 49.4 (100.6) \\
    \hline
  \end{tabular}
  \caption{Weak scaling efficiencies of the M2L, P2P subroutines and
    total algorithm. Efficiencies obtained after adjusting for
    $C_{near}$ and $C_{far}$ are shown in parentheses.}
  \label{tab:weakscaleff}
\end{table}

\subsection{Load Balancing and Adaptivity}
\label{subsec:galaxy}

\begin{figure}[hbt]
  \includegraphics[width=\textwidth]{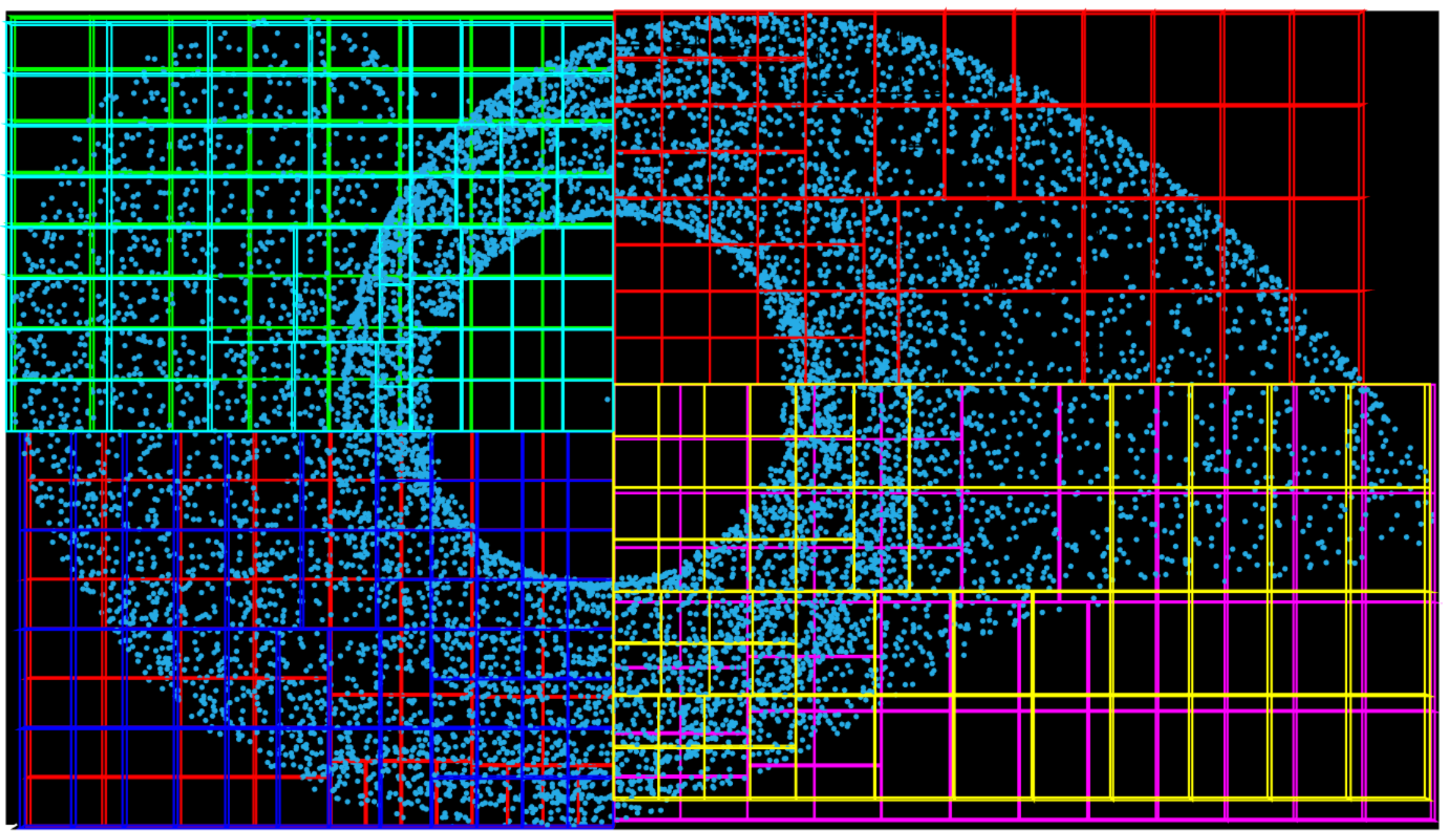}
  \caption{Initial 10,000-point spiral galaxy superimposed on the
    boxes at level 2 generated with $\eta=1$. Colours indicate 8
    different processes.}
  \label{fig:galaxysmallp8}
\end{figure}

As a more complex test case we run a fairly realistic spiral galaxy of
$N$ = 1B $(10^9)$ points using 512 processes and a 5-level tree.  The
distribution was designed to give the model galaxy a multiscale
structure representative of real spiral galaxies. The distribution
forces the balanced tree FMM (with $\eta>0$) to create a tree with
boxes of varying aspect ratios and sizes in close proximity.
We generated the final distribution from a set of 10,000 points in a
rough toroidal shape with inner radius $R_{inner} \approx$ 6000 and
outer radius $R_{outer} \approx$ 18000 parsecs.
This initial set was generated by a Simulink model of two colliding
galaxies \cite{matlabgalaxymodel,toomre:1972}.
Figure~\ref{fig:galaxysmallp8} shows the initial set of 10,000 points
and the boxes on level 2 generated with $\eta=1$. Boxes are coloured
according to which process they belong to ($P$=8 in this image).
The axis of rotation of the galaxy is perpendicular to the page.

The final distribution of $10^9$ points was generated recursively
from the initial set.
For each of the 10,000 points a random distribution of 10 new points
was created within an oblate spheroid centred on that point.
The major axis of each spheroid was $R_0=R_{outer}/20$ and the minor
axis was $r_0=R_0/10$ and parallel to the axis of rotation of the galaxy. For each of the 100,000 points in the new set (the original
points having been deleted) the process was repeated with 
$R_1 = \tfrac 2 3 R_0$ and $r_1 = R_1/10$,
and so on until there were $10^9$ points.
The final set was partitioned into 512 subdomains arranged
as an $8 \times 8 \times 8$ cube of subdomains containing 1,953,125
points each. For simplicity our model only computes the force due to
gravitational attraction, does not include a central black hole,
and does not update the point velocities and positions.

The data was generated and the simulations were run on the Lomonosov-2
supercomputer at the Research Computing Center at Moscow State
University, Russia. This computer was ranked 107(59) in the Top500 in
November 2019 (June 2017) and has similar specifications to Rackham
but 14 cores per node instead of 20.

\begin{figure}[hbt]
  \subfigure[subroutine times]{
    \includegraphics[width=0.8\textwidth]{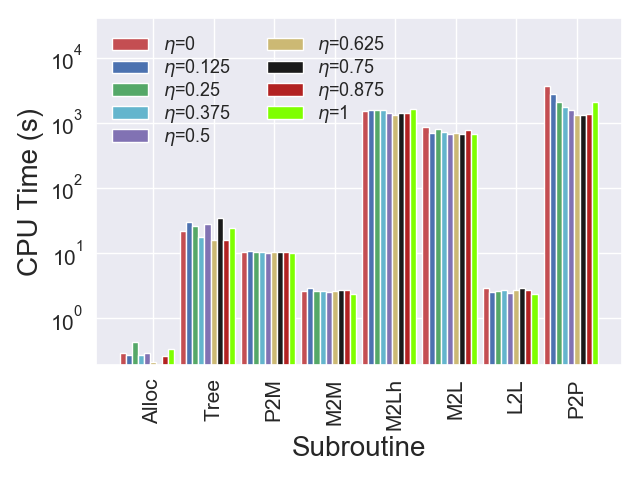}}\\
  \subfigure[total time]{
    \includegraphics[width=0.8\textwidth]{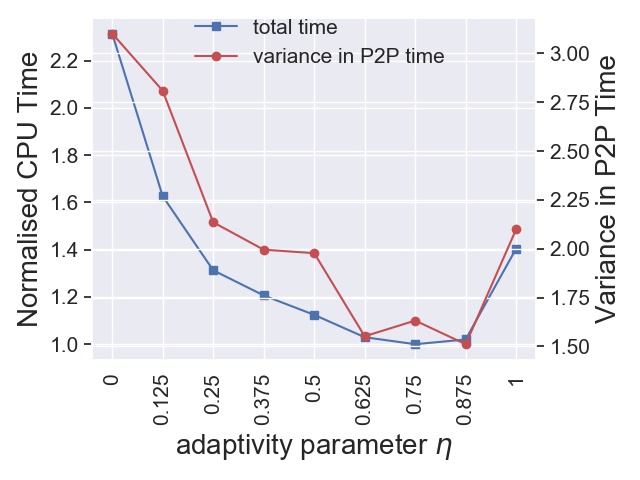}}
  \caption{Sensitivity of galaxy simulation time to $\eta$.
  Simulation used $10^9$ points, 512 processes and a 5 level tree.}
  \label{fig:galaxyl5p512}
\end{figure}

Figure~\ref{fig:galaxyl5p512} (a) shows individual subroutine times.
Note that M2L times are shown separately for the serial M2L shift,
labelled M2L (Algorithm \ref{algorithm:FMMd} lines 17-19), and the
parallel M2L shift, labelled M2Lh (lines 1-16 in Algorithm
\ref{algorithm:FMMd} lines 1-16).  Figure~\ref{fig:galaxyl5p512} (b)
shows the normalised total times (blue) for the galaxy evaluation with
varying $\eta \in [0,1]$.  $\eta$ mostly affects the performance of
P2P in this test.  The value that minimises total time is $\eta =
0.75$ which is 29\% faster than $\eta=1$ (fully adaptive FMM) and 57\%
faster than $\eta=0$ (non-adaptive FMM).
Figure~\ref{fig:galaxyl5p512} (b) also plots the variance in P2P time
(orange), calculated as the difference between maximum and minimum P2P
times divided by the mean P2P time across all processes.  This
quantity is indicative of load balancing.  The lowest variance in P2P
time coincides with the minimum total time, suggesting that load
balancing within P2P is important for efficient execution by reducing
latency.  Optimising over the adaptivity parameter $\eta$ is an
effective way to improve load balancing. Note, however, that even at
$\eta=0.75$ the maximum P2P time is 1.5 times longer than the minimum
so there is room for further improvement.

\section{Conclusion}
\label{sec:conclusions}

A 3D FMM algorithm with distributed parallelism employing balanced
tree adaptivity has been described and tested.  The balanced tree FMM
is highly transparent with several key parameters exposed to user
control. It is also relatively straightforward to implement in
parallel and our code is freely available under a liberal licence at
\url{www.stenglib.org}. We have demonstrated that the performance
response of the balanced tree FMM in parameter space is well-defined
or convex along individual parameters and thus is suited for automated
optimisation.

Strong and weak scaling tests were performed on up to 512 processes
using a modern multicore cluster. Strong scaling efficiency of 23.2\%
was obtained. The low efficiency is primarily due to
the P2P (direct evaluation) subroutine. With a large number of
processes (i.e.~$P \ge 64$) and a relatively small number of points
(10 million), the communication overhead became non-negligible
compared to the computational cost of evaluating interactions, despite
the use of non-blocking MPI send and receive commands. This is a
hardware-related problem rather than an algorithmic one. Better strong
scaling efficiency could be expected if more points were used because
the latency would be hidden by higher computational cost. In this test,
the number of points was limited by the single-core memory capacity.

Weak scaling efficiency of between 46\% to 63\% was obtained depending
on the number of levels. Note that it is not possible to obtain 100\%
weak scaling efficiency due to a geometric effect that results in
increasing connectivity as the domain is divided into more subdomains.
By accounting for this, we estimated an adjusted weak scaling
efficiency of between 73\% and 100\% depending on the number of levels.
Absolute computation times in the weak scaling tests were sensitive to
the number of levels (i.e.~to the number of points per leaf box) and
displayed a convex response. With 1 million points per process the
optimal values were 5 levels or equivalently 244 points per leaf box.

We investigated the effect of the adaptivity parameter $\eta$ on
performance. $\eta=0$ corresponds to non-adaptive FMM (purely
geometric box splitting), balancing the cost of the M2L shift.
Conversely $\eta=1$ corresponds to fully adaptive FMM (purely median
splitting), balancing the cost of the P2P direct evaluation.
Intermediate values, $0<\eta<1$, are thus a blend of both.  The time
taken to simulate gravitational attraction in a large-scale
non-homogeneous point distribution had a convex response to $\eta$,
meaning that it can be readily optimised. The optimal value
$\eta=0.75$ reduced computation time by a factor of more than 2
relative to the non-adaptive case ($\eta=0$) and by a factor of 1.4
relative to the purely adaptive case ($\eta=1$). The speedup is mainly
attributed to reduced variation in the number of points per leaf box
resulting in better load balancing and lower latency in the P2P direct
interaction. Again we note that the convexity means that automatic
optimisation could be beneficial, as was demonstrated previously in 2D
\cite{holm2014}.

Scaling performance of the Tree subroutine (halo matrix construction)
was poor although it made a negligible contribution to total time in
the range of $P$ tested. The poor scaling was the result of
receiving data in process number order, causing
significant latency at large $P$.
The dominant M2L and P2P subroutines used a more sophisticated approach
whereby data was processed in the order in which the non-blocking
receive commands were completed.
P2P exhibited poor scaling when the number of levels was large owing to
having a large number of leaf boxes. This caused the cluster's network
to be flooded with a lot of small messages. Previous work
\cite{bull:ppam} testing the code on the Lomonosov-2 machine suggested
that the scalability limits are system-specific. More software
engineering measures could be taken to improve the code performance.

For best performance in the most important context of time-dependent
problems, continuous autotuning of parameters can be expected to be
highly advantageous, as indicated by the earlier experiences in 2D
\cite{holm2014}. An implementation issue is here to maintain and
locally update the distributed FMM tree at a near optimal cost. This
is also were the transparency of the implementation becomes a definite
advantage.


\section*{Acknowledgment}

This work was financially supported by the Swedish Research Council
within the UPMARC Linnaeus center of Excellence (S.~Engblom), by the
Swedish strategic research programme eSSENCE (J.~Bull, S.~Engblom),
and partly supported by the The Swedish Foundation for international
Cooperation in Research and Higher Education (STINT) Initiation Grant
IB2016-6543, entitled `Large scale complex numerical simulations on
large scale complex computer facilities - identifying performance and
scalability issues', 2016--2017. We are grateful to the Research
Computing Center at the Moscow State University, Russia for granting
us access to Lomonosov-2 where the galaxy simulations were
performed. The scaling tests were performed on the Rackham cluster at
UPPMAX provided by the Swedish National Infrastructure for Computing
(SNIC) and we are grateful for their support. We are also grateful to
Masters students Anders G{\"a}rden{\"a}s and David Ryman of Uppsala
University for their work on testing early versions of the
\texttt{daFMM3D} code.


\newcommand{\doi}[1]{\href{http://dx.doi.org/#1}{doi:#1}}
\newcommand{\available}[1]{Available at \url{#1}}
\newcommand{\availablet}[2]{Available at \href{#1}{#2}}

\bibliographystyle{abbrvnat}
\bibliography{daFMM3D}


\appendix

\section{Algorithms}

Algorithms for the construction of the tree, connectivity matrix and
halo matrix, the downward pass of the balanced FMM tree, and the
direct evaluation are presented here.

\begin{algorithm}
  \caption{Tree, halo and connectivity matrix construction on process $p$.}
  \label{algorithm:FMMh}
  \begin{algorithmic}[1]
    \Statex \textit{Root level ($l$=0)}:
    \State Construct root box $T_p^{(1)}$ and connectivity matrix $C_p^{(1)}$
    \For{$q=0:P-1,q \neq p$}
    \State Non-blocking send centre coords and radius of root box $T_p^{(1)}$ to $q$
    \State Blocking receive centre coords and radius of root box $T_q^{(1)}$ from $q$
    \If{theta\_crit($T_p^{(1)}$,$T_q^{(1)}$)} 
    \State $T_p^{(1)}$, $T_q^{(1)}$ weakly connected
    \Else 
    \State $T_p^{(1)}$, $T_q^{(1)}$ strongly connected
    \EndIf
    \EndFor
    \Statex \textit{Loop over remaining levels}:
    \For{l=1:L}
    \State Adaptive split of parents in $T_p^{(l-1)}$ into 8 children in $T_p^{(l)}$ according to $\eta$ value
    \For{$q=0:P-1,q \neq p$}
    \State Compute no. of connections: \texttt{ncon} $= 64 \times$ strong connections in
    $H_{pq}^{(l-1)}$
    \State Allocate space for lists of \texttt{ncon} connections: $H_{pq}^{(l)}$.\texttt{ibox},\texttt{jbox}
    \State Make \texttt{ilist}, \texttt{jlist}: children
    of strongly connected parents on $p$, $q$ resp.
    \State Non-blocking send centre coords and radii of \texttt{ilist} to $q$
    \State Non-blocking receive centre coords and radii of \texttt{jlist} from $q$
    \EndFor
    \State Construct local connectivity matrix $C_p^{(l)}$ on level $l$
    \For{$q=0:P-1,q \neq p$}
    \State Wait for receive from $q$ to complete
    \State Initialise indices: \texttt{fwd}=0, \texttt{bwd}=\texttt{ncon}
    \For{$i \in$ \texttt{ilist}}
    \For{$j \in$ \texttt{jlist}}
    \If{theta\_crit($i,j$)} 
    \State $T_{p,i}^{(l)}$, $T_{q,j}^{(l)}$ weakly connected
    \State Set $H_{pq}^{(l)}$.\texttt{ibox[bwd]}=$i$, $H_{pq}^{(l)}$.\texttt{jbox[bwd]}=$j$
    \State \texttt{bwd}=\texttt{bwd}-1
    \Else 
    \State $T_{p,i}^{(l)}$, $T_{q,j}^{(l)}$ strongly connected
    \State Set $H_{pq}^{(l)}$.\texttt{ibox[fwd]}=$i$, $H_{pq}^{(l)}$.\texttt{jbox[fwd]}=$j$
    \State \texttt{fwd}=\texttt{fwd}+1
    \EndIf
    \EndFor
    \EndFor
    \EndFor
    \EndFor
\end{algorithmic}
\end{algorithm}

\begin{algorithm}
  \caption{Downward pass of balanced FMM on process $p$.}
  \label{algorithm:FMMd}
\begin{algorithmic}[1]
    \Statex \textit{Parallel M2L}:
		\For{$l=1:L-1$}
		\For{$q=0:P-1, q \neq p$}
        \State From $H_{pq}^{(l)}$ get boxes in $p$ weakly connected to boxes in $q$
        \State Start non-blocking send of box centre coords and coeffs from $p$ to $q$
        \State Start non-blocking receive box centre coords and coeffs from $q$ to $p$
        \EndFor
        \State Set $W$ = $\{q:q \in [0:P-1]\ |\ q$ contains boxes weakly connected to $p \}$
        \While{$W$ is not empty}
		\For{$q \in$ $W$}
		\If{non-blocking receive from $q$ has completed}
        \State Asymmetric M2L shift from boxes in $q$ to boxes in $p$
        \State Remove $q$ from set: $W$ = $W$ - $q$
        \EndIf
        \EndFor
        \EndWhile
        \EndFor
    \Statex \textit{M2L}:
		\For{$l=2:L-1$}
        \State Symmetric M2L shift between boxes in $p$
        \EndFor
    \Statex \textit{L2L}:
		\For{$l=1:L-1$}
        \State L2L shift from boxes in $p$ to their children
        \EndFor
\end{algorithmic}
\end{algorithm}

\begin{algorithm}
  \caption{Evaluation stage of balanced FMM on level $L$ on process
    $p$.}
  \label{algorithm:FMMe}
\begin{algorithmic}[1]
    \Statex \textit{Parallel P2P setup}:
	\For{$q=0:P-1,q \neq p$}
    \State Get list of boxes in $p$ strongly connected to boxes in $q$
    \State Get point coordinates and masses in owned boxes
    \State Start non-blocking send of box indices, coordinates and masses to $q$
    \State Start non-blocking receive of box indices, coordinates and masses from $q$
    \EndFor
    \Statex \textit{L2P (far-field contribution)}:
    \State Evaluate local expansion at points in leaf boxes in $p$
    \Statex \textit{P2P (near-field contribution)}:
    \State Symmetric direct interactions between points in leaf boxes
    in $p$
    \Statex \textit{Parallel P2P evaluation}:
    \State Set $S$ = $\{q:q \in [0:P-1]\ |\ q$ contains boxes strongly connected to $p \}$
    \While{$S$ is not empty}
	\For{$q \in$ $S$}
	\If{non-blocking receive from $q$ has completed}
    \State Asymmetric direct interaction with points in leaf boxes in $q$
    \State Remove $q$ from set: $S$ = $S$ - $q$
    \EndIf
    \EndFor
    \EndWhile
\end{algorithmic}
\end{algorithm}

\end{document}